\documentclass[12pt]{article}
\usepackage[left=2cm,right=2cm,top=1.5cm,bottom=1.5cm]{geometry}
\usepackage{amsfonts}
\usepackage{amssymb}
\usepackage{amsmath}
\usepackage{eufrak}
\usepackage[dvips]{graphicx}
\usepackage{latexsym}
\usepackage{color}

\newcommand{\ind}{\makebox[1em]{\raisebox{-.5ex}[0ex][0ex]{\makebox[0em]%
{$\smile$}}\raisebox{.4ex}[0ex][0ex]{\makebox[-.02em]{$|$}}}}

\newcommand{\C}{{\EuFrak C}}

\newcommand{\bi}{\begin{itemize}}
\newcommand{\ei}{\end{itemize}}

\newtheorem{theorem}{Theorem}[section]
\newtheorem{lemma}[theorem]{Lemma}
\newtheorem{fact}[theorem]{Fact}
\newtheorem{claim}{Claim}
\newtheorem{corollary}[theorem]{Corollary}
\newtheorem{proposition}[theorem]{Proposition}
\newtheorem{definition}[theorem]{Definition}
\newtheorem{remark}[theorem]{Remark}

\newcommand{\bd}[1]{\textbf{#1}}

\newcommand{\be}{\begin{enumerate}}
\newcommand{\ee}{\end{enumerate}}

\newcommand{\ssk}{\smallskip}
\newcommand{\bs}{\bigskip}
\newcommand{\ms}{\medskip}

\newcommand{\cdd}{\cdots}

\newcommand{\bimp}{\leftrightarrow}
\newcommand{\Bimp}{\Leftrightarrow}

\newcommand{\sbb}{\subseteq}

\newcommand{\ya}{\vDash}

\newcommand{\LL}{\mathcal{L}}

\newcommand{\al}{\alpha}
\newcommand{\bt}{\beta}

\newcommand{\og}{\omega}

\newcommand{\acl}{\operatorname{acl}}

\newcommand{\vph}{\varphi}

\newcommand{\scl}{\operatorname{scl}}

\newcommand{\dk}[1]{\emph{#1}}
\newcommand{\set}[1]{ \{ \hspace{0.012in} #1 \hspace{0.012in} \} }
\newcommand{\defi}{\mathrel{\mathop:}=}

\newcommand{\nin}{\noindent}

\def\Ind#1#2{#1\setbox0=\hbox{$#1x$}\kern\wd0\hbox to 0pt{\hss$#1\mid$\hss}
\lower.9\ht0\hbox to 0pt{\hss$#1\smile$\hss}\kern\wd0}

\def\ind{\mathop{\mathpalette\Ind{}}}

\def\Notind#1#2{#1\setbox0=\hbox{$#1x$}\kern\wd0\hbox to 0pt{\mathchardef
\nn=12854\hss$#1\nn$\kern1.4\wd0\hss}\hbox to
0pt{\hss$#1\mid$\hss}\lower.9\ht0 \hbox to
0pt{\hss$#1\smile$\hss}\kern\wd0}

\title{A preservation theorem for theories without the tree property of the first kind}

\author{Jan Dobrowolski\footnote{The first author was supported by NCN Grant no. 2015/19/D/ST1/01174 and 
by Samsung Science Technology Foundation
under Project Number
SSTF-BA1301-03.} and Hyeungjoon Kim}
\date{}

\begin{document}

\setcounter{section}{-1}
\maketitle
\begin{abstract}
 We prove that the NTP$_1$ property of a geometric theory $T$ is inherited by theories of lovely pairs and $H$-structures associated
 to $T$. We also provide a class of examples of nonsimple geometric NTP$_1$ theories.
\end{abstract}

\section{Introduction}

One theme of research in model theory is to  inquire whether some well-known properties are preserved under a certain unary 
predicate expansions of a given structure. One of the motivations for this is that positive theorems of this kind often allows us to obtain interesting and
complicated-looking theories which still satisfy some strong tameness conditions.

The study of expansions by unary predicates reaches back to the paper of Poizat on beautiful pairs \cite{10}, and has been
developed in various directions ever since. Remarkable papers on this subject include for example \cite{11},
where the stability condition is examined, and \cite{12}, which is dedicated to studying expansions in simple 
theories (which generalize the
stable ones).

The well-known equivalence TP $\Bimp$ TP$_1 \vee $ TP$_2$ due to Shelah \cite{5} (where TP denotes the tree property while TP$_1$ and TP$_2$ denote the tree properties of the first and second kind, respectively),  suggests two natural generalizations of simple theories, namely NTP$_1$ theories and NTP$_2$ theories (i.e.,  theories \emph{without}  TP$_1$ and TP$_2$, respectively).   So far, NTP$_1$ and NTP$_2$ theories have been studied much less
 extensively than the simple ones (i.e. theories without TP). However, recently some interesting results on these theories began 
 to appear, notably \cite{3} and \cite{4}. In particular, natural examples of non-simple NTP$_1$ thoeries
 were provided in \cite{3}, namely:
 $\omega$-free PAC fields, linear spaces with a generic bilinear form and a class of theories
 obtained by the ``pfc'' construction.
 
The study of expansions in the NTP$_2$ context was undertaken in \cite{2}, where it was shown that the  NTP$_2$ property 
is preserved under  ``dense and codense'' unary predicate expansions
of geometric structures where  the unary predicate is assumed to  define either an algebraically independent subset or an elementary 
substructure. In the present paper, we prove that the NTP$_1$ property  is 
also preserved under such expansions. One of the 
main ingredients in our proof is the recently proved fact (due to Chernikov and Ramsey \cite{3}) that the TP$_1$ property can, in any TP$_1$ theory, always be witnessed by some formula in a \emph{single} free variable. We also prove (in Section \ref{pfc}) 
that an NTP$_1$ nonsimple geometric theory can be obtained from any Fra{\"\i}ss{\'e} limit which has a simple theory
by applying some constructions from \cite{3}. This yields a large class of nonsimple (so also not NTP$_2$) theories satisfying the assumptions 
of our main theorem.

Our paper is organized as follows. In Section \ref{prelim}, we review  some essential facts about dense codense predicate expansions 
from \cite{2}. In Section \ref{artemoverview}, we state slightly modified versions of some results from \cite{3} concerning SOP$_2$ (equivalently TP$_ 1$). In Section \ref{dense}, we prove our main result, namely, that NTP$_1$ is preserved 
under the unary predicate expansions
defined in Section \ref{prelim}.  In Section \ref{pfc}, we  show that the ``pfc'' construction from \cite{3} preserves a certain
strengthening of geometricity, and conclude from this that a class of geometric nonsimple NTP$_1$ structures can be obtained via  imaginary cover and ``pfc'' operations.

\bs

\section{Dense codense predicate  expansion}\label{prelim}

In this section, we review some basic facts about the dense codense predicate expansions. 

\ssk

Recall that  a theory is called  \dk{geometric} if (1) it eliminates the quantifier $\exists^{\infty}$, and (2) the algebraic closure  
satisfies the exchange property. Examples of geometric theories include all SU-rank 1 theories 
(in particular, strongly minimal theories). For nonsimple NTP$_1$ examples see Section \ref{pfc}.

Throughout, unless stated otherwise, variables may have an arbitrary length. The symbol $\ind$ denotes the algebraic independence relation.

\begin{definition} \label{hstruc}
Let $T$ be a geometric complete theory in a language $\LL$, and let $\LL_H\defi \LL\cup \set{H}$ be the extended language obtained by adding a new unary predicate symbol $H$. For any model $M\ya T$, let $(M, H(M))$ denote an expansion of $M$ to the extended language $\LL_H$, where $H(M)\defi \set{x\in M \mid H(x)}$.
\be
\item We say $(M, H(M))$ is \dk{dense codense}  if every non-algebraic 1-type in $\LL$ over any finite dimensional 
subset $A\sbb M$ has realizations both in $H(M)$ and in $M\setminus \acl_T(A\cup H(M))$.
\item A dense codense expansion $(M, H(M))$ is called a \dk{lovely pair} (resp.\ \dk{H-structure}) if $H(M)$ happens
to be an elementary substructure (resp.\ algebraically independent subset) of $M$.
\ee
\end{definition}

\noindent \bd{Fact/Definition} \cite{7, 8} Let $T$ be any geometric complete theory. Then all of its lovely pairs have 
the same theory, i.e., they are elementarily equivalent to one another. The same holds for $H$-structures. 
We let $T_P$ and $T^{ind}$ denote the common complete theories of the lovely pairs and $H$-structures,
respectively, associated with $T$. By $T^*$, we  shall mean either $T_P$ or $T^{ind}$.

\bs

For the remainder of this section, we shall work inside some fixed, sufficiently saturated model $(M, H(M))\ya T^*$ unless stated otherwise. When $x$ is a tuple of variables, $H(x)$ shall mean the conjunction $H(x_1)\wedge\cdd \wedge H(x_n)$ where $x_i$'s are the variables occurring in $x$. When $A$ is a subset of $M$, $H(A)$ denotes the set $\set{x\in A\mid H(x)}$.

\begin{definition}
For any subset $B\sbb M$, we define  \[\scl(B)\defi \acl_T(B\cup H(M))
\]
which is called the \dk{small closure of $B$}. If $A$ is any subset of $\scl(B)$, we shall say that $A$ is \dk{$B$-small}.
\end{definition}

\begin{definition}
A subset $A\sbb M$ is said to be \emph{$H$-independent} if  $A\ind_{H(A)}H(M)$.
\end{definition}

The following two facts will be important tools in the proof of our main result.


\begin{fact}[\cite 2]\label{Hrest}
For any $\LL_H$-formula $\vph(x, a)$ where $a$ is $H$-independent, there exists some $\LL$-formula $\psi(x, a)$ such that
\[ \vph(x, a)\wedge H(x) \ \bimp\ \psi(x, a)\wedge H(x).
\]
\end{fact}

\begin{fact}[\cite 2]\label{symdif}
For any $\LL_H$-formula $\vph(x, a)$ where $x$ is a single variable and $a$ is $H$-independent, there exists some $\LL$-formula $\psi(x, a)$ such that the symmetric difference $\vph(x, a)\triangle \psi(x, a)$ defines an $a$-small set. 
\end{fact}

We will also use the following observations:

\begin{fact}[\cite 2]\label{indep}
\be
\item For any finite tuple $c$, there exists some finite tuple $h$ in $H(M)$ such that $c\ind_h H(M)$.

\item For any $H$-independent tuple $c$ and any tuple $h$ in $H(M)$, $ch$ is $H$-independent.
\ee
\end{fact}

\ms

We end this section by remarking that all the results and their proofs in this paper may be carried over to many-sorted contexts. However, for the sake of simplifying our arguments,  we shall assume that our theory $T$ is one-sorted throughout the paper.

\section{Overview of some results on SOP$_2$ from \cite{3}}\label{artemoverview}
\label{definition}

In this section, we state some results about SOP$_2$ from \cite{3} in  slightly modified (`localized') versions which  we will need later. But first, let us quickly review some basic terminologies. We consider the language $L_0\defi \{\triangleleft,<_{lex},\wedge\}$ where  $\triangleleft$ and $<_{lex}$ are binary relation symbols and $\wedge$ is a binary function symbol. Then any set $\al^{<\bt}$ (where $\al$ and $\bt$ are ordinals) admits a natural $\LL_0$-structure
whereby  $\triangleleft$ is interpreted as the prefix partial order, $<_{lex}$ as the lexicographic
order and $\wedge$ as the infimum function (with respect to the prefix order).  We will use the following   `localized' version of SOP$_2$:
\begin{definition}\label{sop2}
A formula $\phi(x;y)$ is said to have SOP$_2$ inside  a type $q(x)$ if there are tuples  $(a_\eta)_{\eta\in 2^{<\omega}}$ satisfying the following two properties:
\be
\item For every $\xi\in 2^\omega$, the set $q(x)\cup \set{\phi(x,a_{\xi_|n}):n<\omega}$ is consistent;
\item  For every pair of $\triangleleft$-incomparable elements $\eta, \nu\in 2^{<\og}$, the formula $\phi(x; a_\eta)\wedge \phi(x; a_\nu)$  is inconsistent.
\ee
And a theory has SOP$_2$ inside of $q(x)$ if some formula   has it inside $q(x)$.
\end{definition}

\nin
(The original, non-localized definition of SOP$_2$ is obtained by setting $q=\emptyset$.)

\ms
By compactness, we easily get:
\begin{remark}\label{rem_nonalg}
If $\phi(x;y)$ has SOP$_2$ inside  a type $q(x)$ witnessed by  $(a_\eta)_{\eta\in 2^{<\omega}}$,
then for every $\xi\in 2^\omega$ the type $q(x)\cup \set{\phi(x,a_{\xi_|n}):n<\omega}$ is nonalgebraic.
\end{remark}

Now, let us recall the notion of  \dk{modeling property} on strongly indiscernible trees,
which we will use repeatedly in the paper.

\begin{definition}
We say that a tree $(a_\eta)_{\eta\in S}$ of compatible tuples of elements of a model $M$ is strongly indiscernible
over a set $C\subseteq M$, if 
$$qftp_{L_0}(\eta_0,\dots,\eta_{n-1})=qftp_{L_0}(\nu_0,\dots,\nu_{n-1})$$ implies 
$tp(a_{\eta_0},\dots,a_{\eta_{n-1}}/C)=tp(a_{\nu_0},\dots,a_{\nu_{n-1}}/C)$
for all $n<\omega$ and all tuples \\ $(\eta_0,\dots,\eta_{n-1})$,$(\nu_0,\dots,\nu_{n-1})$ of elements of $S$.
\end{definition}
The following fact comes from \cite{13}.
\begin{fact}\label{based}
Let $\C$ be a monster model of a complete theory. Then for any tree of parameters
$(a_{\eta})_{\eta\in \omega^{<\omega}}$ from $\C$ there is a strongly indiscernible tree $(b_{\eta})_{\eta\in \omega^{<\omega}}$ based on the tree
$(a_{\eta})_{\eta\in \omega^{<\omega}}$,
which means that for every $\eta_0,\dots,\eta_{n-1}\in \omega^{<\omega}$ there are $\mu_0,\dots,\mu_{n-1}\in \omega^{<\omega}$ such that
$qftp_{L_0}(\eta_0,\dots,\eta_{n-1})=qftp_{L_0}(\mu_0,\dots,\mu_{n-1})$ and $tp(b_{\eta_0},\dots,b_{\eta_{n-1}})=tp(a_{\mu_0},\dots,a_{\mu_{n-1}})$.
\end{fact}

Notice that if $q$ is over $\emptyset$, then the 
 consistency condition in Definition \ref{sop2} is preserved under tree modeling. Hence, inside of such a $q$,
SOP$_2$ is always witnessed by a strongly indiscernible tree of parameters. 

\begin{remark}\label{inf}
With the notation from the above definition, the set $q(x)\cup \{ \phi(x,a_{\xi_|n}):n<\omega \}$ has infinitely many realizations
for any $\xi\in 2^\omega$.
\end{remark}

By a very slight modification of the proof of Lemma 4.6 for \cite{3}, we get:

\begin{fact}\label{collapse}
Suppose $(a_\eta)_{\eta\in 2^{<\omega}}$ is a tree strongly indiscernible over $C$ such that $(a_{0^{\alpha}}:0<\alpha<\omega)$
is indiscernible over $cC$. Let 
$$p(y;z)=tp(c;(a_{0\frown 0^\gamma)_{\gamma<\omega}})/C),$$
and let $q(y)$ be a type over $\emptyset$ contained in $tp(c)$. Then, if 
$$q(y)\cup p(y;(a_{0\frown 0^{\gamma}})_{\gamma<\omega})\cup p(y;(a_{1\frown 0^{\gamma}})_{\gamma<\omega}) $$
is inconsistent, then $T$ has SOP$_2$ inside of $q$.
\end{fact}
{\em Proof.}
By naming parameters we can assume that $C=\emptyset$. Suppose the
type $q(y)\cup p(y;(a_{0\frown 0^{\gamma}})_{\gamma<\omega})\cup p(y;(a_{0\frown 1^{\gamma}})_{\gamma<\omega})$ is not
consistent. By compactness and indiscernibility there is a formula $\psi\in p$ such that 
$$q(y)\cup \{\psi(y,a_0,\dots,a_{0\frown 0^{n-1}}),\psi(y,a_1,\dots,a_{1\frown 0^{n-1}})\}$$ is inconsistent.
Then as in \cite{3}, the $n$-fold elongation (see Definition 2.6 from \cite{3}) of $(a_\eta)_{\eta\in 2^{<\omega}}$ 
witnesses that $\psi$ has SOP$_2$ inside of $q$.
\hfill $\square$\\

Using the above fact and modifying the proof of Theorem 4.8 from \cite{3} in the same manner as above (i.e., replacing
any set of formulas related to a consistency condition by its union with an appropriate type over $\emptyset$), we obtain:

\begin{fact}\label{cher}
Suppose a theory $T$ has SOP$_2$ inside of some type $q(x_0,\dots,x_{n-1})=\bigcup_{i<n}q_i(x_i)$. Then, for some $i<n$,
$T$ has SOP$_2$ inside of $q_i(x_i)$.
\end{fact}

\section{The main result}\label{dense}
The aim of this section is to prove our main result, i.e. Theorem \ref{main}. 
First, we give a characterization of SOP$_2$ thoeries that we will need later.

\begin{proposition}\label{prop}
A theory $T$ has $SOP_2$  if there is a formula $\phi(x,y)$ and a strongly indiscernible
tree $(a_\eta)_{\eta\in 2^{<\omega}}$,
such that the set $\{ \phi(x,a_{0^n}):n<\omega \}$ has infinitely many realizations, and
the formula $\phi(x,a_0)\wedge\phi(x,a_1)$ 
has finitely many realizations.
\end{proposition}
{\em Proof.}
Let $\phi(x,y)$ and $(a_\eta)_{\eta\in 2^{<\omega}}$ be as above.
For a set $B\subseteq 2^{<\omega}$, we put 
$$A_B:=\bigcap_{\eta\in B}\phi(\C,a_\eta),$$ and
for a tuple $b=(\eta_0,\dots,\eta_{n-1})\in (2^{<\omega})^n$,
we put $$A_b:=\bigcap_{i<n}\phi(\C,a_{\eta_i}).$$
\begin{claim}
We can assume that \[(\forall n,m\in \omega\backslash \{0\})(A_{0,1}=A_{0^n,1^m}) \tag{*}. \]
\end{claim}
{\em Proof of Claim 1.}
Since the set $A_{\{0^k,1^k:k=1,2,3,\dots\}}$ is finite, by compactness, it is equal to
$D:=A_{\{0^k,1^k:k=1,2,3,\dots\,K\}}$ for some $K<\omega$. So for any positive, pairwise distinct
$k_1,...,k_{2K}<\omega$, the set $A_{0^{k_1},\dots,0^{k_K},1^{k_{K+1}},\dots,1^{k_{2K}}}$ is contained
in $D$, but by indiscernibility it has the same (finite) cardinality as $D$, so it is equal to $D$.
Replacing $\phi(x,y)$ by $\psi(x,y_0,\dots,y_{K-1})=\bigwedge_{i<K}\phi(x,y_i)$ and
each $a_{\epsilon_0,\dots,\epsilon_{m-1}}$ by $a_{\epsilon_0^K,\dots,\epsilon_{m-1}^K}$
we obtain a tree with the desired properties.
\hfill $\square$\\
So we will assume that $(*)$ holds.

\begin{claim}
We can  assume (in addition to $(*)$) that 
\[ A_{\{0^n1:n<\omega\}}=\emptyset \tag{**} . \]
\end{claim}
{\em Proof of Claim 2.}
Since the set $A_{\{0^n1:n<\omega\}}$ is finite,  it is equal to
$D:=A_{\{0^n1:n=0,1,2,\dots\,K-1\}}$ for some $K<\omega$.
Then $D=A_{0^{k_0}1,\dots,0^{k_{K-1}}1}$ for any pairwise distinct $k_i$'s. Since the tree $(a_{0^K\frown \eta})_{\eta\in
2^{<\omega}}$ is strongly indiscernible over $\{a_{0^k1}:k=0,1,2,\dots\,K-1\}$,
we can replace $\phi(x,y)$ by $\phi(x,y)\wedge\neg\psi(x,y_0,\dots,y_{K-1})$, where 
$\psi(x,a_1,\dots,a_{0^{k-1}1})$ defines $D$, and $a_\eta$ by 
$(a_{0^K\frown \eta},a_1,\dots,a_{0^{K-1}1})$, guaranteeing $(**)$ while preserving $(*)$ and
 strong indiscernibility of the tree.
\hfill $\square$\\
Let us assume that $(*)$ and $(**)$ hold.
Then there is a maximal $n$ such that $A_{\{ 1,01,\dots,0^{n-1}1, 0^{n},0^{n+1},\dots  \}}$ is
 nonempty.
 Put $c=(a_1,\dots,a_{0^{n-1}1})$, $\psi(x,y_0,\dots,y_{n-1})=\bigwedge_{i<K}\phi(x,y_i)$ and $\phi'(x,y,y_0,\dots,
 y_{n-1})=\phi(x,y)\wedge
 \psi(x,y_0,\dots,y_{n-1})$. We claim that the strongly indiscernible tree 
 $(b_\eta)_{\eta\in{2^{<\omega}}}$,
 where $b_\eta:=a_{0^n\frown\eta }c$, 
 witnesses
 SOP$_2$ of $\phi'(x,y)$. Indeed, 
 by the choice of $n$ (and by strong indiscernibility) all paths are consistent.
 Moreover, by maximality of $n$, the set $\{\phi'(x,b_1c)\}\cup \{\phi'(x,b_{0^k}c):k=1,2,\dots\}$
 is inconsistent, but by $(*)$ (and by the strong indiscernibility of the tree $(a_\eta)$) 
 this set is equivalent to the formula $\phi'(x,b_0c)\wedge\phi'(x,b_1c)$, so we get that the 
 latter formula is inconsistent, and we are done.
 (Note that, \emph{a posteriori}, by Remark \ref{rem_nonalg}, the $n$ chosen above must be 
 equal to zero, i.e. the formula $\phi(x,a_0)\wedge\phi(x,a_1)$ is already inconsistent if we 
 assume $(*)$ and $(**)$.)
\hfill $\square$\\

For the remainder of this section, we will work inside a sufficiently saturated model $(M, H(M))\ya T^*$.

\begin{lemma}
Let $\phi(x,y)$ be any $\LL_H$-formula witnessing SOP$_2$. Then for some dummy variables z, the formula $\phi(x,yz)$
witnesses SOP$_2$ with some strongly indiscernible tree consisting of $H$-independent tuples.
\end{lemma}
{\em Proof.}
Let $(a_{\eta})_{\eta\in 2^{<\omega}}$ be a strongly indiscernible tree witnessing that $\phi(x,y)$ is SOP$_2$.
By Fact \ref{indep}, choose a finite tuple $h_{\emptyset}$ of elements of $H(M)$ such that $a_\emptyset h_\emptyset$ is $H$-independent.
For any $\eta\in 2^{<\omega}$ let $h_\eta$ be a conjugate of $h_\emptyset$ under an automorphism (in the sense of $T_H$)
sending $a_\emptyset$ to $h_\emptyset$. Then any indiscernible tree based on $(a_\eta h_\eta)_{\eta\in 2^{<\omega}}$ will satisfy the conclusion.
\hfill $\square$\\

We will need one more preparatory lemma.

\begin{lemma}\label{on_H}
If there is some $\LL_H$-formula $\phi(x,y)$ such that $\phi(x,y)\wedge H(x)$ witnesses SOP$_2$,
then $T$ has SOP$_2$.
\end{lemma}
{\em Proof.}
By \ref{cher} (applied to $\phi(x,y)$ and $q(x):=H(x)$), we can assume that $x$ is a single variable.
Let $(a_\eta)_{\eta\in 2^{<\omega}}$ be a strongly indiscernible tree witnessing SOP$_2$ of $\phi(x,y)\wedge H(x)$ such that
$(a_\eta)$ is $H$-independent. By Fact \ref{Hrest}, there is an $\LL$-formula $\psi(x,y)$ agreeing with $\phi$ on $H$.
For any $\eta\in 2^{<\omega}$, the formula $\psi(x,a_{\eta\frown 0})\wedge\psi(x,a_{\eta\frown 1})$ is algebraic, since otherwise, by 
the density of $H$, it would be realized inside of $H$, a contradiction.
So, by Proposition \ref{prop}, $\psi(x,y)$ has SOP$_2$.
\hfill $\square$\\


In the final proof we will use one more characterization of TP$_1$ property, which was proved in \cite{1}.
First, we remind the definition of $k$-TP$_1$ from there:

\begin{definition}
A formula $\psi(x,y)$ has $k$-TP$_1$ if there are tuples $c_\beta$, $\beta\in \omega^{<\omega}$  such that
for each $\beta\in \omega^\omega$ the set $\{\psi(x,c_{\beta_{|m}}):m<\omega\}$ is consistent,
and for any pairwise incomparable
elements $\beta_0,\dots,\beta_{k-1}\in \omega^{<\omega}$ the set $\{\psi(x,c_{\beta_i}:i<k\}$ is inconsistent.
\end{definition}

\begin{fact}\label{ktp}
Suppose that an $\LL$ formula $\phi(x,y)$ and a tree $(a_\eta)_{\eta\in \omega^{<\omega}}$ witness $k$-TP$_1$
for some $k\geq 2$. Then for some $d<\omega$, the $\LL$-formula 
$\psi(x,y_0,\dots,y_d)=\phi(x,y_0)\wedge\dots\wedge \phi(x,y_d)$
witnesses $2$-TP$_1$.
\end{fact}

Now we are in a position to prove the main result.

\begin{theorem}\label{main}
If $T^*$ has SOP$_2$, then so does $T$.
\end{theorem}
{\em Proof.}
Assume $T^*$ has SOP$_2$ witnessed by an $\LL_H$-formula $\phi(x,y)$, where $x$ is a single variable (we can assume that by \ref{cher})
and a strongly indiscernible tree $A=\{(a_\eta)_{\eta\in 2^{<\omega}}\}$, where each $a_\eta$ is $H$-independent. \
\newline
Case 1: No realization of $\wedge_i \phi(x,a_{0^i})$ is in scl(A).\
By Fact \ref{symdif}, let $\psi(x,y)$ be a formula such that for each $\eta$, $\phi(x,a_\eta)\triangle\psi(x,a_\eta)$ defines
an $a_\eta$-small set. Then for any $\eta$,
$\psi(x,a_{\eta\frown 0})\wedge \psi(x,a_{\eta\frown 1})$ has finitely many realizations, 
since otherwise, by the co-density condition, it would have a realization outside of $scl(A)$, so realizing
the formula $\phi(x,a_{\eta\frown 0})\wedge \phi(x,a_{\eta\frown 1})$.
Also, every realization of $\wedge_i \phi(x,a_{0^i})$ is a realization of $\wedge_i \psi(x,a_{0^i})$, 
so we are done by Proposition \ref{prop}.\
\newline
Case 2: There is some $b\in scl(A)$ satisfying $\wedge_i \phi(x,a_{0^i})$.\
So $b$ realizes some algebraic formula $\theta(x,c,h)$, where $c$ and $h$ are tuples of elements of $A$ and $H$, respectively.
We can assume that for any $c'$ and $h'$ the formula $\theta(x,c',h')$ has at most $k$ realizations, where $k<\omega$ is fixed.
Choose $N<\omega$ such that $c$ is contained in $\{a_\eta:\eta\in 2^{<N}\}$.  Put $d_\eta:=a_{0^N\frown\eta}$. Then $\phi(x,y)$
together with the tree $(d_\eta)_{\eta\in 2^{<\omega}}$, which is strongly indiscernible over $c$, still
witnesses SOP$_2$.
Put $$\mu(z,c,y):=H(z)\wedge \exists x(\theta(x,c,z)\wedge\phi(x,y)).$$
Then, since $\wedge_n \mu(z,c,d_{0^n})$ is realized by $h$ (this is is witnessed by substituting $b$ for $x$),
we get by the indiscernibility of $(d_\eta)_\eta$ over $c$ that
$\wedge_n \mu(z,c,d_{\xi_{|n}})$ is consistent for any $\xi\in 2^\omega$.
Also, for any pairwise incomparable $\eta_1,\dots,\eta_n\in 2^{<\omega}$,
the set $\{\mu(z,c,d_{\eta_i}):i\leq n\}$ is $k+1$-inconsistent. Hence, by compactness,
$\mu(z,x,y)$ has $(k+1)-TP_1$. It follows from Fact \ref{ktp} that some $\LL_H$ formula of the form $H(z)\wedge\nu(z)$
has $TP_1$, so also SOP$_2$. We conclude by Lemma \ref{on_H}.

\hfill $\square$\\

\section{Examples of geometric nonsimple NTP$_1$ theories.}\label{pfc}
\bs

We start by outlining the ``pfc'' construction from Subsection 6.3 of \cite{3}. For the reader's convenience, we repeat the definitions 
used there.
\begin{definition}
Suppose $K$ is a class of finite structures. We say that $K$ has the Strong Amalgamation Property (SAP) if given 
$A,B,C\in K$ and embeddings $e:A\to B$ and $f:A\to C$ there is $D\in K$ and embeddings $g:B\to D$ and 
$h:C\to D$ such that\\
1) $ge=hf$ and\\
2) $im(g)\cap im(h)=im(ge)$ (and hence $=im(hf)$ as well).
 \end{definition}
We will say that a theory is SAP if it has a countable ultrahomogeneous model whose age is SAP.
The following criterion comes from \cite{9}.
\begin{fact}\label{algebraicity}
Suppose K is the age of a countable structure M. Then the following are equivalent:\\
1) $K$ has SAP\\
2) $M$ has no algebraicity
\end{fact}

Let $K$ denote an SAP Fra{\"\i}ss{\'e} class in a finite relational language $\LL=(R_i:i<k)$, 
where each $R_i$ has arity $n_i$.
Denote by $T$ the theory of the Fra{\"\i}ss{\'e} limit of the class $K$.
Then $\LL_{pfc}$ is defined to be a two-sorted language, with the sorts denoted by $O$ and $P$, and relation symbols
$R^i_x(x,y_1,y_2,\dots,y_{n_i})$, where $x$ is a variable of the sort $P$ and $y_i$'s are variables of the sort $O$.
Given an $\LL_{pfc}$-structure $M=(A,B)$ and $b\in B$, the $\LL$-structure associated to $b$ in $M$, denoted $A_b$,
is defined to be the $\LL$-structure interpreted in $M$ with domain $A$ and each $R_i$ interpreted as $R^i_b(A)$.
Put $$K_{pfc}=\{ M=(A,B)\in  Mod(L_{pfc}): |M|<\omega, (\forall b\in B) (\exists D\in K) (A_b\simeq D)\}.$$
\begin{fact}[\cite{3}]
$K_{pfc}$ is a Fra{\"\i}ss{\'e} class satisfying SAP.
\end{fact}
Thanks to the above fact, there is a unique countable ultrahomogeneous $\LL_{pfc}$-structure with age $K_{pfc}$.
Let $T_{pfc}$ denote its theory. Then $T_{pfc}$ has quantifier elimination.
Let us recall two facts from \cite{3} that will be crucial for us.
\begin{fact}\label{A_b}
Suppose $(A,B)\models T_{pfc}$. Then, for all $b\in B$, $A_b\models T$.
\end{fact}
\begin{fact}\label{nonsimp}
Suppose $T$ is a simple theory which is the theory of a Fra{\"\i}ss{\'e} limit of a SAP Fra{\"\i}ss{\'e} class $K$. Then $T_{pfc}$
is $NSOP_1$. Moreover, if the $D$-rank of $T$ is at least 2, then $T_{pfc}$ is not simple.
\end{fact}
Now, we aim to prove that the ``pfc'' construction applied to a geometric theory satisfying condition $acl(A)=A$
for any $A$, gives a geometric theory.

Let $N=(A,B)$ be a monster model of $T_{pfc}$.
\begin{lemma}\label{l1}
For any $A_0\subseteq A$ and $B_0\subseteq B$ we have that $acl(A_0B_0)\cap B=B_0$.
\end{lemma}
{\em Proof.}
Clearly we can assume that both $A_0$ and $B_0$ are finite.
Put $C=A_0B_0$, take any $b\in B\backslash B_0$ and fix any natural number $n$. We will show that the orbit of
$b$ over $C$ has at least $n$ elements. To see this, consider a finite $\LL_{pfc}$-superstructure $E=(A_0,D)$ of 
$(A_0,B_0)\subseteq N$, where
$D=B_0\cup\{b,d_1,\dots,d_n\}$ with $d_j$'s being pairwise distinct elements not belonging to 
$\{b\}\cup B_0$, and for each $i,j$, $R^i_{d_j}(A_0)$ is equal to $R^i_{b}(A_0)$ in the sense of $N$. Then clearly 
$E\in K_{pfc}$, so we can assume $E=(A_0,D)$ is a substructure of $N$ (by changing $d_i$'s appropriately).
Then, for every $j$, we have that $qftp(d_j/C)=qftp(b/C)$, so we are done due to the quantifier elemination in $T_{pfc}$.
\hfill $\square$\\

\begin{lemma}\label{l2}
Suppose $\phi_j(x,a_j,b_j)$ for $j=1,\dots,n$ are non-algebraic $\LL_{pfc}$-formulas in a single variable of the first sort, where
each $a_j$ is a tuple of elements of $A$ and $b_j$ are pairwise distinct elements of $B$. 
Then the conjunction $\phi:=\wedge_j \phi_j$ is non-algebraic.
\end{lemma}
{\em Proof.}
Suppose for a contradiction that $\phi$ is algebraic, and denote by $A_0$ the finite set $\phi(A)$.
We consider a finite $\LL_{pcf}$-structure $E$ with universe $(C,D)$, where $C=A_0\cup \{c\}$, 
$c\notin A_0$, $D=\{b_1,\dots,b_n\}$,
and interpretation of symbols of $\LL_{pcf}$ given as follows. For any $j\leq n$, by the non-algebraicity of 
$\phi_i$
we choose its realization $c_j\in A\backslash A_0$. Now, let $f_j:C\to A_0\cup\{c_j\}$ be the bijection
whose restriction to $A_0$ is the identity. We interpret every $R^i_{b_j}$ in $E$ as $f_j^{-1}[R^i(A_0\cup\{c_j\})]$.
Then $E$ belongs to  $K_{pfc}$ so it embeds in $N$ via some function $g$, and $g(c)\in N\backslash A_0$
is an realization of $\phi$. This is a contradiction to the choice of $A_0$, so the lemma is proved.
\hfill $\square$\\

\begin{corollary}\label{preservation}
If $T$ is geometric and satisfies condition $acl(A)=A$ for any $A$, then $T_{pfc}$ is geometric satisfying the same condition.
\end{corollary}
{\em Proof.}
First, we show the definability of infinity. If $\phi(x,y)$ is a formula with $x$ being a single variable of
the sort $P$ (and $y$ of any lenght), then by Lemma \ref{l1}, for any $c$, if $\phi(x,c)$ is algebraic then
each of its realizations belongs to $c$, so $\phi(x,c)$ is algebraic iff 
it has at most $|y|$ realizations. Now, if $x$ is a variable of the sort $O$, then any formula in variable $x$ can be 
presented in the form $\phi(x,y_{j,l},z_{j,l})_{j,l}:=\vee_l\wedge_j \phi_{j,l}(x,y_{j,l},z_{j,l})$,
where each $z_{j,l}$ is a single variable
of the sort $P$, and $y_{j,l}$ are tuples of variables of the sort $O$ (we can obtain such a presentation
since atomic formulas in $\LL_{pfc}$ can involve only a single variable from the sort $P$).
Then $\phi(x,a_{j,l},b_{j,l})_{j,l}$ is algebraic if and only if for each $l_0$ the formula
$\wedge_j \phi_{j,l_0}(x,a_{j,l_0},b_{j,l_0})$ is algebraic. But by Lemma \ref{l2} this holds iff there are $j_1,\dots,j_s$
such that $b_{j_1,l_0}=b_{j_2,l_0}=\dots=b_{j_s,l_0}$ and $\wedge_{t\leq s} \phi_{j,l_0}(x,a_{j_t,l_0},b_{j_t,l_0})$ is algebraic. 
By Fact \ref{A_b} the latter is a definable condition on $a_{j,l_0},b_{j,l_0}$, so we obtain the definability of infinity.

As to the condition $acl(A)=A$, by Lemma \ref{l1} it is enough to check that for any finite $A_0\subset A$ and $B_0\subset B$
we have that $acl(A_0B_0)\cap A=A_0$. Consider any $a\in acl(A_0B_0)\cap A$. Then there are formulas $\phi_j(x,a_j,b_j)$ as in the 
statement of Lemma \ref{l2}, such that the conjunction $\phi(x):=\wedge_j \phi_j(x,a_j,b_j)$ is algebraic and satisfied by $a$.
By Lemma \ref{l2}, for some $j$ the formula $\phi_j(x,a_0,b_j)$ is algebraic. By Fact $\ref{A_b}$ and the assumptions on
$T$ this implies that $a$ in $acl(a_j)=a_j$ in the sense of the structure $A_{b_j}$, so $a\in A_0$.

\hfill $\square$\\

\begin{remark}
 By a similar argument we can show, assuming only that $T$ is geometric, that in $T_{pfc}$ we have definability of
 infinity, and the following weaker form of exchange principle:
 $$a\in acl(Ab)\backslash acl(A) \implies b\in acl(Aa).$$
 for any parameter set $A$ and $a$ belonging to the same sort as $b$. However, this condition seems not to be sufficient 
 to prove a generalization of Theorem \ref{main} by our methods.
\end{remark}

Let us recall one more operation from \cite{3} which we need to obtain examples of nonsimple NTP$_1$ theories.
Given an $\LL$-structure $M$, the imaginary cover $\tilde{M}$ of $M$ is defined to be the structure in language $\LL'$ obtained from $\LL$
by adding a binary relation symbol $E$, constructed by replacing each element of $M$ with an infinite $E$-class and interpreting
the symbols of $\LL$ in the natural way. By Remark 6.19 from \cite{3} we have:

\begin{fact}\label{D>1}
If $T=Th(M)$ is simple and SAP, then $\tilde{T}:=Th(\tilde{M})$ is simple of $D$ rank at least 2 and SAP.
\end{fact}
Let us notice the following.
\begin{remark}\label{tilde}
The theory $\tilde{T}$ has definability of infinity and for any $A$ in a model of $\tilde{T}$ we have that $acl(A)=A$.
\end{remark}
{\em Proof.}
The second clause is obvious, and for the first one notice that any atomic formula $\phi(x,a)$ in $\tilde{T}$ 
has either infinitely many or at most one realization. Hence, by the quantifier elimination we get definability of infinity for
any formula.
\hfill $\square$\\

Now we obtain the final corollary which yields a class of examples of nonsimple NTP$_1$ geometric theories.

\begin{corollary}
If $T$ is any SAP simple theory then $(\tilde{T})_{pfc}$ is a geometric nonsimple NTP$_1$ theory.
\end{corollary}
{\em Proof.}
By Remark \ref{tilde}, $\tilde{T}$ is geometric and satisfies the condition $acl(A)=A$ for any $A$, so by Corollary \ref{preservation}
the same is true about $(\tilde{T})_{pfc}$. Moreover, $(\tilde{T})_{pfc}$ is NTP$_1$ and nonsimple by 
Facts \ref{nonsimp} and \ref{D>1}.
\hfill $\square$\\

\noindent
{\bf Addresses:}\\
Jan Dobrowolski\\
Department of Mathematics, Yonsei University\\
50 Yonsei-ro, Seodaemun-gu, Seoul 120-749, South Korea\\\
{\bf E-mail:} dobrowol@math.uni.wroc.pl \\[2mm]
Hyeungjoon Kim\\
Department of Mathematics, Yonsei University\\
50 Yonsei-ro, Seodaemun-gu, Seoul 120-749, South Korea\\
{\bf E-mail:} joon9754@yonsei.ac.kr\\

\end{document}